# A roadmap to the (vanishing of the) Euler characteristic

CLARA LÖH AND GEORGE RAPTIS

The Euler characteristic $\chi(X) \in \mathbb{Z}$ of a finite CW-complex $X$ is a classical and fundamental homotopy invariant of $X$. It is the unique numerical homotopy invariant of finite CW-complexes that is additive and its value at the one-point space is 1. For an oriented closed connected smooth $n$-manifold $M$, the Poincaré dual of the Euler characteristic $\chi(M) \in H_0(M;\mathbb{Z})$ is the Euler class $e(M) \in H^n(M;\mathbb{Z})$ of the tangent bundle of $M$; the vanishing of the characteristic class $e(M)$ (or $\chi(M)$) is equivalent to the existence of a nowhere vanishing vector field on $M$. Moreover, the classical Poincaré–Hopf theorem states that the Euler characteristic is equal to the sum of the indices of a vector field at its isolated zeroes. By Poincaré duality, the Euler characteristic of an odd-dimensional closed manifold vanishes.

The Euler characteristic and its vanishing appear prominently in various contexts and in connection with different kinds of invariants and structures, especially, in the case of aspherical manifolds. More specifically, the properties of the simplicial volume or the $L^2$-Betti numbers of $M$, the existence of affine structures on $M$ or flat connections on its tangent bundle, as well as the existence of Riemannian metrics on $M$ with specific constraints on the sectional curvature, all bear (established or conjectural) relations to the vanishing (or not) of the Euler characteristic. Not only are these relations interesting and significant, but also the various ways they are themselves interconnected and support each other in different degrees of intensity are intricate and intriguing, too.

In order to illustrate this, let us choose as starting point a well-known conjecture attributed to Chern [25, 24]: *the Euler characteristic of a closed affine manifold vanishes.* This is known to hold for complete affine manifolds [25] and for special affine manifolds [24]. (A conjecture of Markus states that a closed affine manifold is complete if and only if it is special affine; see [15].) Complete affine manifolds are aspherical and their universal covering space is $\mathbb{R}^n$ – the last property fails for general closed aspherical $n$-manifolds [11]. The tangent bundle of a general affine manifold $M$ supports a flat connection. Therefore, as a consequence of the generalized Milnor–Wood inequality [23, 6, 14], the Euler characteristic of $M$ vanishes if the simplicial volume of $M$ vanishes.

This suggests the conjecture that already *the simplicial volume of a (complete) affine manifold vanishes* [5], which was shown for complete affine manifolds whose holonomy contains a pure translation [5]. In addition, the last conjecture and the Chern conjecture are interestingly related with a question of Gromov [17][29]: *does the vanishing of the simplicial volume of an aspherical*

*Date*: September 8, 2023. This work was supported by the CRC 1085 Higher Invariants (Universität Regensburg, funded by the DFG)., We are grateful to Grigori Avramidi and Bruno Klingler for helpful discussions.

*manifold imply the vanishing of its Euler characteristic?* Note that this implication is false for non-aspherical manifolds (e.g. $S^2$). Gromov's question is also equivalent to the boundedness of the Euler class of an aspherical manifold [29] – which brings us back to properties of the tangent bundle. It becomes tempting to think that the existence of a flat connection on the tangent bundle (not necessarily torsion-free) might already imply the vanishing of the Euler characteristic, but this turns out to be false in general [36].

On the other hand, if $M$ admits a flat Riemannian metric, then $M$ is aspherical and also both the simplicial volume and the Euler characteristic of $M$ vanish. More generally, if $M$ admits a Riemannian metric of non-negative sectional curvature, then its simplicial volume vanishes [9, 16] and a conjecture attributed to Hopf predicts that its Euler characteristic is non-negative (see [31]). Of course such manifolds are not aspherical in general. Instead, by the classical Cartan–Hadamard theorem, $M$ is aspherical if it admits a Riemannian metric of non-positive curvature and, in addition, the simplicial volume of $M$ is positive if $M$ is actually negatively curved [21].

From another viewpoint, the Euler characteristic can also be computed from the $L^2$-Betti numbers.

A well-known conjecture attributed to Singer states (see [31]): *the $L^2$-Betti numbers of an aspherical manifold vanish away from the middle dimension.* This problem evolved in several stages: the starting point was Singer's work on $L^2$-harmonic forms [35] and the continuation was due to Dodziuk [13] and Anderson [1]. (The analogue of the Singer conjecture over finite fields has recently been disproved [3].) Combined with Gromov's question, the Singer conjecture asks whether *the $L^2$-Betti numbers of an aspherical manifold vanish if its simplicial volume vanishes.* Using Gromov's *Main Inequality* [16], the latter condition follows from the vanishing of the minimal volume. In this case, that is, when $M$ is aspherical with vanishing minimal volume, then its $L^2$-Betti numbers indeed vanish [33], and therefore so does its Euler characteristic, as a consequence of the $L^2$-Euler–Poincaré formula [31].

One is left to marvel at the complex pattern drawn by these implications and to feel intrigued by their intuitive or formal relations. The diagram below aims to be a useful condensed guide to these results and conjectures, summarized non-linearly in the form of a network of interconnected (established or conjectural) implications. We recommend the interested reader to consult the references provided in the diagram for further reading.

## CONTENTS

## 1. Glossary

**1.1. Standing assumption.** In the diagram, $M$ is an oriented closed connected (smooth) manifold of dimension $n$.

**1.2. Notation.**

**1.** $\chi(M)$: the *Euler characteristic* of $M$

**2.** $e(M)$: the *Euler class* of the tangent bundle of $M$; its $\ell^\infty$-norm $\big\|e(M)\big\|_\infty$ is measured via bounded cohomology [16, 14]

**3.** $\mathrm{minvol}(M)$: the *minimal volume* of $M$ [16, 4]

**4.** $\|M\|$: the *simplicial volume* of $M$ [16, 27]

**5.** $\|M\|_{\mathcal{F},\mathbb{Z}}$: the *integral foliated simplicial volume* of $M$ [34, 18, 28]

**6.** $\widetilde{M}$: the *universal covering* of $M$

**1.3. Terminology.**

**7** (affine). An *affine* manifold is a smooth manifold with an atlas whose transition maps are affine transformations. Equivalently, an affine manifold is a smooth manifold with a flat torsion-free connection on its tangent bundle [2].

**8** (almost abelian). A group $\Gamma$ is *almost abelian* if there exists a finite normal subgroup $N \lhd \Gamma$ such that $\Gamma/N$ contains an abelian subgroup of finite index. In particular, every almost abelian group is amenable [32].

**9** (bounded). A cohomology class is *bounded* if it lies in the image of the comparison map from bounded cohomology [16, 22] to ordinary cohomology.

**10** (complete affine). A smooth affine manifold $M$ is *complete affine* if the associated developing map $\widetilde{M} \to \mathbb{R}^n$ is a homeomorphism.

**11** (flat). A Riemannian metric on a smooth manifold is *flat* if it has constant sectional curvature 0. On complete connected smooth manifolds, a Riemannian metric is flat if and only if the Riemannian universal covering space is isometric to a Euclidean space.

**12** (non-negative curvature). A Riemannian metric on a smooth manifold has *non-negative curvature* if its sectional curvature is everywhere $\geq 0$.

**13** (non-positive curvature). A Riemannian metric on a smooth manifold has *non-positive curvature* if its sectional curvature is everywhere $\leq 0$.

**14** (negative curvature). A Riemannian metric on a smooth manifold has *negative curvature* if its sectional curvature is everywhere $< 0$.

**15** (pure translation). An affine manifold *has a pure translation* if its holonomy contains a pure translation.

**16** (small amenable category). A manifold $M$ has *small amenable category* if it admits an open cover by at most $\dim(M)$ amenable subsets [16, 7] (or, equivalently, an open cover by amenable subsets with multiplicity at most $\dim(M)$).

**17** (special affine)**.** A smooth affine manifold $M$ is *special affine* if its holonomy lies in the special affine group [24].

**18** (tangentially flat)**.** A smooth manifold is *tangentially flat* if its tangent bundle admits a flat connection.

### 1.4. **Arrows.**

**19** (arrow types)**.**

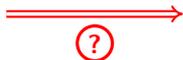 conjectural implication

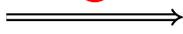 straightforward implication / by definition

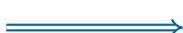 non-trivial implication / by theorem

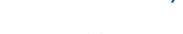 this implication does *not* hold in general

⊖ this implication also holds without "asphericity"

C. Löh. Fakultät für Mathematik, Universität Regensburg, Regensburg, Germany, clara.loeh@ur.de, https://loeh.app.ur.de

G. Raptis. Fakultät für Mathematik, Universität Regensburg, Regensburg, Germany, georgios.raptis@ur.de, https://graptismath.net/

**Boxes**

- the Euler class of $M$ is bounded
- $M$ is aspherical
- the $L^2$-Betti numbers of $M$ vanish away from the middle dimension
- $\widetilde{M} \cong \mathbb{R}^n$
- $M$ is complete affine
- $M$ is complete affine and has a pure translation
- $M$ is aspherical and $\|M\| = 0$
- $M$ is aspherical and $\|M\|_{\mathcal{F},\mathbb{Z}} = 0$
- $M$ is special affine
- $\|M\| = 0$
- $M$ is aspherical and $\mathrm{minvol}(M) = 0$
- $M$ is aspherical and has small amenable category
- $M$ is affine
- $\chi(M) = 0$
- all $L^2$-Betti numbers of $M$ vanish
- $M$ admits a metric of non-positive curvature
- if $\dim M = 2 \cdot d$, then $(-1)^d \cdot \chi(M) \geq 0$
- $\|e(M)\|_\infty \leq 1/2^n$ and $|\chi(M)| \leq \|M\|/2^n$
- $M$ is tangentially flat
- $M$ is tangentially flat and $\|M\| = 0$
- $M$ admits a flat metric
- $M$ admits a metric of negative curvature
- $\|M\| > 0$
- $\pi_1(M)$ is amenable
- $\pi_1(M)$ is almost abelian
- $M$ admits a metric of non-negative curvature
- $\chi(M) \geq 0$

**Edge labels**

- reformulation of Gromov's question [29]
- Singer conjecture [31, Ch 11][35, 13, 1]
- [1], 12]
- Cartan–Hadamard theorem
- [5]
- [34]
- Mahler conjecture
- [25]
- Gromov's question [17, p. 232][29]
- vanishing theorem [16, 22]
- [36]
- [34]
- [16, 4]
- main inequality
- Gromov [17, p. 232]
- [16, Sec 3.4]
- [34]
- [5]
- [24]
- $\|M\| = 0$
- Chern–Gauß–Bonnet [16, 14]
- Chern conjecture [25, 24, 19]
- [18, § 5.34][33]
- $L^2$-Euler–Poincaré formula [31, Theorem 1.35]
- $L^2$-Euler–Poincaré formula [31, Theorem 1.35]
- Hopf conjecture [20][31, Conj 11.2]
- generalized Milnor–Wood inequality [23, 6][14, Sec 13.3]
- generalized Milnor–Wood inequality
- Bieberbach classification [8]
- [16][26, Cor 6.4][14, Sec 8.14]
- straightening [21]
- mapping theorem [16]
- [32, Ch 0]
- [9]
- Hopf conjecture [20][31, Conj 11.2]